\documentclass[11pt,amssym,twoside]{article}
\usepackage{amssymb}
\newtheorem{thm}{Theorem}[section]

\newtheorem{prop}[thm]{Proposition}

\newtheorem{conj}[thm]{Conjecture}

\newcommand{\Cplx}{\mathbb C}

\begin{document}

\title{Approximation on abelian varieties by its subgroups}%
\author{Arash Rastegar}%


\maketitle
\begin{abstract}
In this paper, we introduce an algebro-geometric formulation for Falting's theorem on diophanitine approximation on abelian varieties
using an improvement of Faltings-Wustholz observation over number fields [Fal-Wus]. In fact, we prove that, for any geometrically irreducible subvariety $E$
of an abelian variety $A$ and any finitely generated subgroup $F$ of $A(\mathbb{C})$ we have an estimate of the form $d_v(E; x) \gg H(x)^{-\delta }$
for the distance $d_v (E; x)$ of a point $x$ in $F$
outside $E$ where $v$ is a place of $K$.
This was proved before, only for $F$ being the set of rational points of $A$ over a number field [Fal].
\end{abstract}
\section*{Introduction}

Historically diophantine approximation has been a strong method for 
proving finiteness results in diophantine geometry. 
In the case of the projective line the fundamental result is Roth's theorem. 
In this paper, we shall generalize this to geometrically irreducible subvarieties of abelian varieties.
We fix a finitely generated field $K$ of characteristic zero and chose a geometrically irreducible subvariety $E\subset \mathbb{P}^n$
defined over $K$ and study the distance $d_v (E; x)$ of a $K$-rational point 
of $\mathbb{P}^n-E$ to $E$ where $v$ is a place of $K$. It is natural to ask for an estimate 
of the type 
$$
d_v(E; x) \gg H(x)^{-\delta }
$$
such that the positive exponent $\delta$ is as small as possible. Here $H(x)$ denotes the absolute 
Well height of $x$. Faltings and Wustholz showed that the above assertion is trivial for subvarieties $E$
which are geometrically irreducible.

\begin{thm}
Let $A$ be an abelian variety defined over a finitely generated subfield $K$ of $\mathbb{C}$. Let
$E$ is a geometrically irreducible subvariety of $A$ defined over $K$ and $F$ be a finitely generated subgroup of 
$A(K)$.                                                                                                                                                                                                                                                                       
Let $v$ be a valuation on $K$ and $H(x)$ a height function on $K$ coming from a choice of projective model for $K$
over the algebraic closure of $\mathbb{Q}$ in $K$.
If $d_v(x,E)$ denotes the $v$-adic distance from $x$ to $E$, and $\delta$ and $c$ are positive constants,
then, there are only finitely many points in $F$ satisfying the following inequality
$$
d_v(x,E)< cH(x)^{-\delta }.
$$
\end{thm}

This implies that Faltings' theorem on Diophantine approximation on abelian varieties is
an algebro-geometric fact, not an arithmetic one. Here is a geometric implication in the case of affine subsets of abelian varieties.

\begin{thm}
Let $A$ be an abelian variety
defined over $\mathbb{C}$ and let $U$ be an open affine subset of $A$. Suppose
$F\subset \mathbb{C}^n$ is a finitely generated subgroup.
Then $F\cap U(\mathbb{C})$ is a finite set.
\end{thm}

This extends Lang's conjecture on integral points of affine subvarieties of abelian varieties.
Lang had the same geometric expectation when he formulated his conjecture that a curve of genus $\geq 2$ in its Jacobian should
intersect any finitely generated subgroup of the Jacobian in a finite set.
He even conjectured that divisible group of any finitely generated subgroup of the Jacobian intersects the embedded curve 
in finitely many points. Getting such a result is beyound our reach. Since we have only access to a result for finitely generated subgroups
defined over a finitely generated field. Therefore we state a conjecture following geometric philosophy of Lang as follows:

\begin{conj}
Let $A$ be an abelian variety
defined over $\mathbb{C}$ and let $U$ be an open affine subset of $A$. Suppose
$F\subset \mathbb{C}^n$ is a finitely generated subgroup. Let $Div(F)$ denote the divisible subgroup of $\mathbb{C}^n$ associated to $F$.
Then $Div(F) \cap U(\mathbb{C})$ is a finite set.
\end{conj}

\section{Diophantine approximation by subgroups of $\mathbb{C}^n$}

This section is devoted to proving theorems whichhas motivated us
to extend Faltins' result. This section is borrowed from [Ras1]. The arguments are along the same lines as
analogous classical results.

Roth's theorem on Diophantine approximation of rational points on
projective line implies a version on projective varieties defined
over number-fields. 

\begin{thm} (Improvement of Roth's theorem on diopphantine approximation)
Fix a finitely generated field of characteristic zero $K$ 
and $\sigma :K\hookrightarrow \Cplx$ a
complex embedding. Let $A$ be an abelian variety
defined over $K$ and let $L$ be an very ample line-bundle on
$A$. Denote the arithmetic height function associated to the
line-bundle $L$ by $h_L$. Suppose $F\subset A(K)$ is a finitely generated subgroup.
Fix a Riemannian metric on $A_{\sigma}(\Cplx)$ and
let $d_{\sigma}$ denote the induced metric on
$A_{\sigma}(\Cplx)$. 
Then, for every $\delta>0$ and every choice
of an algebraic point $\alpha\in A(\bar {K})$ 
and all choices of a
constant $C$, there are only finitely many points
$\omega\in F$ approximating $\alpha$ such that 
$$
d_{\sigma}(\alpha ,\omega)\leq Ce^{-\delta h_L(\omega)}.
$$
\end{thm}

\begin{prop}
With assumptions of the above theorem, suppose for some $\delta_0>0$ 
we have that, for any choice of a constant $C$
and every choice
of an algebraic point $\alpha\in A(\bar {K})$
there are only finitely many points
$\omega\in F$ approximating $\alpha$ in the following manner
$$
d_{\sigma}(\alpha ,\omega)\leq Ce^{-\delta_0 h_L(\omega)}.
$$
Then, for every $\delta>0$ and every choice
of an algebraic point $\alpha\in A(\bar {K})$ and all choices of a
constant $C$, there are only finitely many points
$\omega\in F$ approximating $\alpha$ such that 
$$
d_{\sigma}(\alpha ,\omega)\leq Ce^{-\delta h_L(\omega)}.
$$
\end{prop}
\textbf{Proof (Proposition).} Note that, we have assumed that the above is true for some
$\delta_0>0$ 
without any assumption  
on $\alpha$.
Let $\delta'>0$ be infimum of such $\delta_0>0$.
The subset $F$ is disjoint union of the images of finitely many height-increasing
self-endomorphisms $\phi_i:A(K)\to A(K)$ 
defined over $K$ such that for
all $i$ we have
$$
h_L(\phi_i(f))=mh_L(f)+O(1)
$$
where $m>1$. Take $\phi_i:A(K)\to A(K)$ to be of the form $u \mapsto mu+a_i$
where $a_i$ are representatives of the finite group quotient $F/mF$.

Fix $\epsilon>0$ such that $\epsilon<\delta' <m\epsilon$. 
Suppose that $w_n$ is an infinite sequence of elements in
$F$ such that $\omega_n\to \alpha$ which satisfies the estimate
$$
d_{\sigma}(\alpha ,\omega_n)\leq Ce^{-\epsilon h_L(\omega_n)}.
$$
then infinitely many of them are images of elements of $F$ under
the same $\phi_i$. By going to a subsequence, one can find a
sequence $\omega'_n$ in $F$ and an algebraic point $\alpha'$ in
$A(\bar {K})$ such that $\omega'_n \to \alpha'$ and for a fixed
$\phi_i$ we have $\phi_i(\alpha')=\alpha$ and
$\phi_i(\omega'_n)=\omega_n$ for all $n$. Then
$$
d_{\sigma}(\alpha ,\omega_n)\leq Ce^{-\epsilon h_L(\omega_n)}\leq
C'e^{-\epsilon m_i h_L(\omega'_n)}
$$
for an appropriate constant $C'$. On the other hand,
$$
d_{\sigma}(\alpha' ,\omega'_n)\leq C''d_{\sigma}(\alpha ,\omega_n)
$$
holds for an appropriate constant $C''$ and large $n$ by
injectivity of $d\phi_i^{-1}$ on the tangent space of $\alpha$.
This contradicts our assumption on $\delta'$, because $\delta' <m_i\epsilon$.
$\square$ 
\\
\textbf{Proof (Theorem).}  If we assume that points of $F$ and covering maps are defined over some
number-field, Roth's theorem implies that the assumption of theorem is true for any 
$\delta_0>2$. The same is true for finitely generated fields of characteristic zero
by a result of Lang [Lan] generalizing Roth's theorem.$\square$

\section{Diophantine approximation on abelian varieties}

The version of Roth's theorem in the last section being true, 
we expect the following version of Liouville's theorem on diophantine approximation to hold:

\begin{thm}(Weak version of Vojta conjectures on diophantine approximation)
Fix a finitely generated field of characteristic zero $K$. Let $V$ be a smooth projective algebraic variety
defined over $K$ and let $L$
be an very ample line-bundle on $V$. Denote the arithmetic height
function associated to the line-bundle $L$ by $h_L$. 
Then, there exists a positive constant $\delta_0$ such that for any positive constant $c$
For any geometrically irreducible algebraic subvariety $E$ of $V$ defined over $K$ 
and $d_v(x,E)$ denoting the $v$-adic distance from $x$ to $E$,
there are only finitely many points defined over $K$ in $V(K)$ outside $E$ satisfying the following inequality
$$
d_v(x,E)< cH(x)^{-\delta_0}
$$
except for points in an algebraic variety $V(\delta_0 )$ which is of strictly smaller dimension of $V$.
In case $K$ is trancendental, we have to pick a model for $K$ over algebraic closure of 
$\mathbb{Q}$ in $K$ following Lang [Lan] and [Lan-Ner] to fix a height function.
\end{thm}
\textbf {Proof.} This is a weak form of Vojta conjectures.
In the number field case, this is mentioned in Faltings-Wustholz [Fa-Wu] as a trivial result in case $E$ is geometrically irreducible.
In the case of finitely generated fields of characteristic zero the result is a consequence of theorem I' in seminal work of Lang [Lan].$\square$

Now, the following version of Falting's theorem, can be proved 
using the method of height expansion.

\begin{thm}(Geometric formulation of Faltings' theorem on diophantine approximation on abelian varieties)
Fix a finitely generated field of characteristic zero $K$. Let $A$ be an abelian variety
defined over $K$ and let $L$
be an very ample line-bundle on $A$. Denote the arithmetic height
function associated to the line-bundle $L$ by $h_L$. Suppose
$F\subset A(K)$ is a finitely generated subgroup
Fix any positive constants $\delta$ and $c$. 
For any irreducible algebraic subvariety $E$ of $A$ defined over $K$ 
and $d_v(x,E)$ denoting the $v$-adic distance from $x$ to $E$,
there are only finitely many points defined over $K$ in $F$ outside $E$ satisfying the following inequality
$$
d_v(x,E)< cH(x)^{-\delta}
$$
Again, in case $K$ is trancendental, we have to pick a model for $K$ over algebraic closure of 
$\mathbb{Q}$ in $K$ following Lang [Lan] and [Lan-Ner] to fix a height function.
\end{thm}
\textbf {Proof.}
We have that the above is true for some
$\delta_0>0$ 
without any assumption  
on $\alpha$.
Let $\delta'>0$ be infimum of such $\delta_0>0$.
Note that a finitely generated subgroup 
$F\subset A(K)$ is a subset which is a union of its images under
finitely many height-increasing polynomial finite self-endomorphisms
$\phi_i:A\to A$ defined over $K$ such that for all $i$ we have
$$
h_L(\phi_i(u)) > mh_L(u)+0(1)
$$
where $m>1$. Let $\phi_i$ be of the form $u\mapsto mu+a_i$ where $a_i$ are representatives of $F/mF$ for $m>1$.
Fix $\epsilon>0$ such that $\epsilon<\delta' <m\epsilon$. 
Suppose that $P_n$ is an infinite sequence of elements in
$F$ such that $P_n\to E$ which satisfies the estimate
$$
d_{\sigma}(E ,P_n)\leq Ce^{-\epsilon h_L(P_n)}.
$$
then infinitely many of them are images of elements of $F$ under
the same $\phi_i$. By going to a subsequence, one can find a
sequence $P'_n$ in $F$ and an irreducible component of inverse image of $E$ under $\phi_i$ denoted by $E'$ in
$A(\bar {K})$ such that $P'_n \to E'$ and for a fixed
$\phi_i$ we have $\phi_i(E')=E$ and
$\phi_i(P'_n)=P_n$ for all $n$. Then
$$
d_{\sigma}(E ,P_n)\leq Ce^{-\epsilon h_L(P_n)}\leq
C'e^{-\epsilon m h_L(P'_n)}
$$
for an appropriate constant $C'$. On the other hand,
$$
d_{\sigma}(E' ,P'_n)\leq C''d_{\sigma}(E ,P_n)
$$
holds for an appropriate constant $C''$ and large $n$ by
injectivity of $d\phi_i^{-1}$ on the tangent space of $E$.
This contradicts our assumption on $\delta'$, because $\delta' <m\epsilon$.
So, reduction of the inequalty for some positive $\delta_0$ to arbitrary $\delta >0$ is done in the 
same manner as in our version of Roth's theorem. $\phi_i$ are finite maps and and if $E$ is approximated by infinitely many points
$P_n$, then a subsequence are in the image of some $\phi_i$ and a therefore an infinite subsequence of its inverse images approximate an irreducible
component of inverse image of $E$ under $\phi_i$.
Getting rid of $A(\delta )$ is the result of the fact that
$A(\delta )$ is invariant under $\phi_i$ and therefore a union of abelian subvarieties and $F\cap A(\delta )$ 
is again a translation of a finitely generated subgroup. One can proceed by reducing the problem from $A$ and $E$
to $A(\delta )$ and $E\cap A(\delta )$ and applying induction, until $A(\delta )$ is union of finitely many points.$\square$

The following will be an important implication:

\begin{thm}(Geometric formulation of Lang's conjecture on diophantine approximation on affine subsets of abelian varieties by integral points)
Fix a finitely generated field of characteristic zero $K$. Let $A$ be an abelian variety
defined over $K$ and let $L$
be an very ample line-bundle on $A$. Denote the arithmetic height
function associated to the line-bundle $L$ by $h_L$. Let $U$ be an open affine subset of $A$. Suppose
$F\subset K^n$ is a finitely generated subgroup.
Then $F\cap U(K)$ is a finite set.
\end{thm}
\textbf {Sketch of Proof.} Assume that $l$ is an equation for $E$. The height $H(x)$, for 
$x \in A - E$ an integral point, is essentially the inverse of the product of the 
$v$-norms of $l(x)$, $v$ running through the infinite places of $K$. Our proof actually 
gives that all these are bounded below by multiples of $H(x)^{-\delta }$, so that $H(x)$ 
must be bounded.$\square$ 

\subsection*{acknowledgements}
I have benefited from conversations with M. Hadian, A. Rajaei, P. Sarnak, N. Talebizadeh, 
for which I am thankful. Peter Sarnak particularly gave crucial comments which led 
to the final version of the paper. 
I would also like to thank Sharif University
of Technology for finantial support and Princeton University for warm hospitality.


Sharif University of Technology, e-mail: rastegar@sharif.edu
\\Princeton University, e-mail:
rastegar@princeton.edu

\end{document}